\documentclass{article}
\usepackage[english]{babel}
\usepackage{amsmath,amssymb,stmaryrd,latexsym,theorem}
\usepackage{enumerate}
\usepackage{a4wide}

\newenvironment{proof}{\noindent\textbf{Proof\ }}{\hspace*{\fill}$\Box$\medskip}
\newtheorem{theorem}{Theorem}
\newtheorem{lemma}[theorem]{Lemma}
\newtheorem{corollary}[theorem]{Corollary}
\newtheorem{assumption}{Assumption}
\newtheorem{definition}{Definition}
{\theorembodyfont{\rmfamily}\newtheorem{remark}{Remark}}

\DeclareMathOperator{\Div}{div}
\DeclareMathOperator{\dist}{dist}

\newcommand{\R}{\mathbb{R}}
\newcommand{\Z}{\mathbb{Z}}
\newcommand{\average}[1]{\langle #1 \rangle}
\newcommand{\eps}{\varepsilon}
\newcommand{\F}{\mathcal{F}}
\newcommand{\I}{\mathcal{I}}
\newcommand{\X}{\mathcal{X}}
\newcommand{\T}{\mathcal{T}_{\eps}}
\newcommand{\Tin}{\mathcal{T}_{\eps}^{\rm in}}
\newcommand{\Tbl}{\mathcal{T}_{\eps}^{\rm bl}}
\newcommand{\E}{\mathcal{E}^{\eps}}
\newcommand{\Ein}{\mathcal{E}^{\eps}_{\rm in}}
\newcommand{\Oin}{\Omega_{\eps}^{\rm in}}
\newcommand{\Obl}{\Omega_{\eps}^{\rm bl}}
\newcommand{\jump}[1]{\llbracket #1 \rrbracket}
\newcommand{\tn}[1]{\interleave #1 \interleave_{\eps}}

\begin{document}

\title{Quantitative homogenization of Stokes flow in a perforated domain under low regularity of data and geometry}
\author{
	Alexei Lozinski \\
	Université Marie et Louis Pasteur, CNRS, LmB (UMR 6623), F-25000 Besançon, France
}
\maketitle

\begin{abstract}
	We study the homogenization of the steady Stokes equations with no-slip conditions in a domain perforated by $\varepsilon$-periodic obstacles, whose
	limit is Darcy's law for the homogenized pressure $p^\ast$. We prove the quantitative error estimate with rate $\sqrt\varepsilon$, natural in the presence of a
	boundary layer, only assuming that $f-\nabla p^\ast\in H^1$. Such estimates were previously available only under strong smoothness 	assumptions, such as $f\in C^{1,1/2}$. Our geometrical assumptions are also less stringent than  in the previous analyses. The outer domain and the reference fluid cell need only be Lipschitz, the perforations may deviate from periodicity in a layer of width $O(\varepsilon)$ near the outer boundary, and the solid part may be
	connected (in three dimensions) and attached to the outer wall.  In the absence of
	a boundary layer, the same argument yields the rate $O(\varepsilon)$ assuming only $f\in H^1$.
\end{abstract}

\section{Introduction}\label{sec:intro}

We consider the Stokes problem in the perforated domain $\Omega^\eps:=\Omega\setminus\bar B^\eps$, $\Omega\subset\R^d$, $d=2$ or 3: find $u^\eps:\Omega^\eps\to\R^d$ and $p^\eps:\Omega^\eps\to\R$ solution of
\begin{alignat}{2}
  -\Delta u^\eps+\nabla p^\eps&=f&\quad&\text{in }\Omega^\eps,\label{genP}\\
  \Div u^\eps&=0&&\text{in }\Omega^\eps,\label{genD}\\
  u^\eps&=0&&\text{on }\partial\Omega^\eps,\label{genB}
\end{alignat}
with a given  $f\in L^2(\Omega)^d$. $B^\eps$ represents periodic perforations placed on a Cartesian grid of period $\eps$. To define them, we start from the unit reference cell $Y=(0,1)^d$, introduce an open set $B\subset Y$ (the reference obstacle) and $\F=Y\setminus\bar B$ (the reference fluid domain); the precise geometric hypotheses are collected in Assumptions~\ref{AssPattern} and~\ref{AssOmeps} below. For $\eps>0$ and $k\in\Z^d$ let $Y_k^\eps=\eps(Y+k)$, $B_k^\eps=\eps(B+k)$, $\F_k^\eps=\eps(\F+k)$. A typical choice is
\begin{equation}\label{setI}
  \I=\{k\in\Z^d:\ Y_k^\eps\cap\Omega\neq\varnothing\},\qquad B^\eps=\Big(\bigcup_{k\in\I}\overline{B_k^\eps}\cap\Omega\Big)^\circ,\qquad \Omega^\eps=\Omega\setminus\bar B^\eps .
\end{equation}
The actual definition of $\Omega^\eps$ may differ from \eqref{setI} in a band of width $\sim\eps$ near the outer boundary $\partial\Omega$.

To study the behavior of the solution for small $\eps$, one introduces the cell Stokes problems: for $i=1,\dots,d$, find $w_i:\F\to\R^d$ and $\pi_i:\F\to\R$, $\Z^d$-periodic, such that
\begin{alignat}{2}
  -\Delta w_i+\nabla\pi_i&=e_i&\quad&\text{in }\F,\label{eq:cell}\\
  \Div w_i&=0&&\text{in }\F,\notag\\
  w_i&=0&&\text{on }\partial B,\notag\\
  \int_\F\pi_i&=0. &&\notag
\end{alignat}
The formal two-scale asymptotic expansion, cf.\ e.g.\ \cite{hornung,Sanchez80}, suggests that $(u^\eps,p^\eps)$ is close to $(u^\ast,p^\ast)$ with
\begin{equation}\label{eq:homvel}
  u^\ast=\eps^2w_i\Big(\frac{\cdot}{\eps}\Big)(f_i-\partial_ip^\ast),
\end{equation}
where $p^\ast$ satisfies the Darcy problem on the whole (non-perforated) domain
\begin{align}
  \Div\big(A^\ast(f-\nabla p^\ast)\big)&=0\quad\text{in }\Omega,\label{darcy}\\
  n\cdot A^\ast(f-\nabla p^\ast)&=0\quad\text{on }\partial\Omega.\label{DarcyB}
\end{align}
The constant matrix $A^\ast$ is defined by $A^\ast e_i=\average{w_i}$, $i=1,\dots,d$, where $\average{\cdot}$ stands for the average on $\F$:
\begin{equation}\label{devavg}
  \average{v}=\frac{1}{|\F|}\int_\F v(y)\,dy.
\end{equation}

The homogenization of the Stokes equations was first rigorously studied in \cite{Tartar80}, where the weak $L^2$ convergence of the velocity and the strong $L^2$ convergence of the pressure were established. The strong $L^2$ convergence of the velocity was later proven in \cite{allaire89}, where the case of a connected $B^\eps$ was also considered (as opposed to the case of $B^\eps$ composed of isolated obstacles strictly inside each periodic cell $Y_k^\eps$ in the earlier literature). More recently, convergence rates of the type
\begin{equation}\label{RateAll}
  \frac{1}{\eps^2}\|u^\eps-u^\ast\|_{L^2(\Omega^\eps)}+\frac1\eps|u^\eps-u^\ast|_{H^1(\Omega^\eps)}+\|p^\eps-p^\ast\|_{L^2(\Omega^\eps)}\leqslant C\sqrt\eps
\end{equation}
have been obtained. The first such estimate is available in \cite{Mik96}, with the sub-optimal rate $\sqrt[6]{\eps}$ instead of $\sqrt\eps$, for $d=2$. The convergence with the correct order $\sqrt\eps$ has been proven in \cite{StokesPart2} for $d=2$ and in \cite{shen22,balazi24,balaziPhD} for any $d$. Apart from their methods of proof (boundary layer correctors in \cite{shen22} vs.\ cut-off techniques elsewhere), these articles differ in their assumptions on the perforations: in \cite{StokesPart2} and \cite{shen22}, $B^\eps$ is composed of isolated obstacles, whereas \cite{balazi24,balaziPhD} cover the technically more involved and physically more realistic case of a connected structure $B^\eps$ in dimension $d=3$. All of these papers make rather strong smoothness assumptions on the data $f$ and on the boundary $\partial\Omega^\eps$. In \cite{shen22}, $f$ is assumed to be in $C^{1,1/2}(\bar\Omega)^d$ and the boundary $\partial B$ (hence $\partial B^\eps$) of class $C^{1,\alpha}$. In \cite{balazi24}, $f-\nabla p^\ast$ is assumed to be in $H^2(\Omega)^d\cap C^{1,\alpha}(\bar\Omega)^d$ and $\partial B^\eps$ of class $C^{m+2}$ with $m>d/2$. Obtaining an $O(\sqrt\eps)$ estimate with less regularity on $f$, than $C^{1,1/2}$ is singled out as ``an interesting and challenging problem'' in \cite{shen22}. Lowering the regularity assumptions on $\partial B^\eps$ is arguably even more important for applications in engineering and in numerical analysis, where $B^\eps$ typically has corners and edges.

The purpose of the present article is to establish an estimate of type \eqref{RateAll} under weaker regularity assumptions on $f$ and $\partial B^\eps$. We only assume that $f-\nabla p^\ast\in H^1(\Omega)^d$. We make no smoothness assumptions on $\Omega$ or on the periodic perforations beyond Lipschitz conditions plus some mild restrictions on the behavior of the perforations near $\partial\Omega$ (where they are allowed to deviate from the periodic pattern). In particular, $B^\eps$ can be connected (in dimension $d=3$) and touch the outer boundary $\partial\Omega$. We note in passing that the preceding works, even those allowing for a connected structure, assume that the structure is separated from the wall $\partial\Omega$ by a free space (filled by fluid) of width $\sim\eps$, which is a non-physical situation of a structure floating freely in the fluid. The only exception is Appendix~B of \cite{balaziPhD}, where a connected $B^\eps$ touching $\partial\Omega$ is allowed, but $\Omega$ is restricted to be a box $\prod_{i=1}^d(0,n_i\eps)$ with integers $n_i$.

Under the assumptions of the present article, the solutions $(w_i,\pi_i)$ of the cell problems have only the regularity $H^1(\F)^d\times L^2(\F)$, and one cannot guarantee that $u^\ast$, as defined by \eqref{eq:homvel}, belongs to $H^1(\Omega^\eps)^d$. We circumvent this problem by replacing $f-\nabla p^\ast$ in the definition of $u^\ast$ with its interpolant by continuous piecewise polynomial functions on the mesh naturally associated with the periodic perforations, i.e. the uniform Cartesian mesh of size $\eps$. We thus introduce the regularized homogenized velocity
\begin{equation}\label{reghomvel}
  \hat u^\ast=\eps^2w_i\Big(\frac{\cdot}{\eps}\Big)\,I^0_\eps(f_i-\partial_ip^\ast),
\end{equation}
where $I^0_\eps$ is a Cl\'ement-type interpolation \cite{clement} onto the $\mathbb Q_1$ finite element space on the $\eps$-cells $Y_k^\eps$. $I^0_\eps$ produces continuous functions  vanishing on $\partial\Omega$ and, if necessary, in a band of width $\sim\eps$ near $\partial\Omega$. By construction, $\hat u^\ast\in H^1_0(\Omega^\eps)^d$, so that the estimate
\begin{equation}\label{RateAlluhat}
  \frac{1}{\eps^2}\|u^\eps-\hat u^\ast\|_{L^2(\Omega^\eps)}+\frac1\eps|u^\eps-\hat u^\ast|_{H^1(\Omega^\eps)}+\|p^\eps-p^\ast\|_{L^2(\Omega^\eps)}\leqslant C\sqrt\eps
\end{equation}
makes sense. This regularization by finite elements is similar to a device introduced in
the theory of operator error estimates for elliptic homogenization
\cite{birman04,zhikov05,ZP05}. There, the classical first-order approximation
$u_0+\eps\chi(\cdot/\eps)\nabla u_0$ is replaced by
$u_0+\eps\chi(\cdot/\eps)S_\eps\nabla u_0$, where $S_\eps$ is the Steklov average over cubes of size $\eps$. The purpose is the same as here: to make sense of the product $\chi(\cdot/\eps)S_\eps\nabla u_0$ under minimal regularity assumptions. The estimate $O(\sqrt\eps)$ in the presence of boundary layers is also proved by these methods in \cite{ZP05,PS12,suslina13}, see also \cite{shen2018}. Using
$\mathbb Q_1$ interpolant in place of $S_\eps$ has certain advantages in the present setting: the cut-off near $\partial\Omega$ is built into the interpolant, which
vanishes outside the interior cells. With $S_\eps$, on the contrary, one needs an extension of $u_0$ beyond $\Omega$ and a separate cut-off function, as in \cite{shen2018}.

Besides the differences in the regularity settings and geometrical assumptions, our contribution differs from the papers on periodic homogenization of Stokes equations on perforated domains mentioned above in its  method of proof. It borrows ideas from the finite element analysis, introducing an inf-sup condition between velocity and pressure with respect to a mesh-dependent pressure norm, proved by a local argument in the spirit of the macroelement technique \cite{Stenberg84,KS92}. This norm includes the jumps of the cell means of the pressure across the faces of the $\eps$-mesh, which is particularly reminiscent of the stabilized $\mathbb Q_1$--$\mathbb P_0$ finite element methods \cite{HF87,SK90,KS92}, cf. also the textbook \cite{ESW14}.

\section{Assumptions and main result}\label{secAssmp}

\begin{assumption}[periodic pattern]\label{AssPattern}
	$B\subset Y=(0,1)^d$ is open and non-empty, $\F=Y\setminus\bar B$ is connected and Lipschitz. Moreover, for $j=1,\dots,d$, the fluid part of two cells	adjacent in the direction $e_j$, i.e.
	$	\F^{(j)}=\Big(\overline{\F}\cup\big(\overline{\F}+e_j\big)\Big)^\circ $,
	is connected and Lipschitz.
\end{assumption}
To formulate the second assumption, we recall the notations $Y_k^\eps=\eps(Y+k)$, $\F_k^\eps=\eps(\F+k)$, for $k\in\Z^d$, and introduce the definition of a special class of domains, in which one can construct the inverse of operator $\Div$, cf. Lemma \ref{lemBog} below:
\begin{definition}[John domains] \label{deftnqohn}
	A bounded domain \(\omega\subset\R^d\) is a John domain if there exist a point \(x_0\in\omega\) (the John center) and a constant \(J>0\) such that, for every \(x\in\omega\), there is a rectifiable curve
	$
	\gamma:[0,L]\longrightarrow\omega
	$
	parametrized by the arc length and satisfying 
	\(\gamma(0)=x\), 
	\(\gamma(L)=x_0\), 
	and 
	\(
		\operatorname{dist}(\gamma(t),\partial\omega)
		\geqslant {t}/{J}.
	\) for every \(t\in[0,L]\).
\end{definition}  
The class of John domains is bigger than that of Lipschitz, but it excludes certain geometrical features, like external cusps or narrow channels connecting two parts of the domain (assuming a moderate value for the John constant $J$).
 
\begin{assumption}[perforated domain and mesh]\label{AssOmeps}
$\Omega\subset\R^d$, $d=2$ or $3$, is a bounded connected Lipschitz domain and $\Omega^\eps\subset\Omega$ is open. Moreover, $\bar\Omega$ is covered by a mesh $\T$ of closed cells with pairwise disjoint interiors, $\bar\Omega=\bigcup_{T\in\T}T$, divided into two disjoint families $\T=\Tin\cup\Tbl$ such that:
\begin{enumerate}[(a)]
\item (interior cells) $\Tin=\{\overline{Y_k^\eps}:k\in\I_0\}$ for a set $\I_0\subset\Z^d$, and $\Omega^\eps\cap Y_k^\eps=\F_k^\eps$ for all $k\in\I_0$;
\item (boundary layer) every $T\in\Tbl$ satisfies $T\subset\bar\Omega$ and $\dist(x,\partial\Omega)\leqslant c_0\eps$ for all $x\in T$;
\item (shape of the layer cells) every $T\in\Tbl$ is of diameter
$\leqslant c_2\eps$;  $\eps^{-1}T$ is Lipschitz with constants independent of $\eps$; the fluid part of any such cell, i.e. $T^\circ\cap\Omega^\eps$, is a John domain of constant $J$ and of has measure $\geqslant c_1\eps^d$; moreover, for any other cell $T'\in\T$ adjacent to $T$, i.e. such that the intersection $T\cap T'$ is of positive $(d-1)$-measure, the fluid part $(T\cup T')^\circ\cap\Omega^\eps$ is also a John domain of constant $J$. There are at most $N_\mathrm{i}$ such cells $T'$ for any $T\in\Tbl$.
\end{enumerate}
The parameters $c_0,c_1,c_2,J,N_\mathrm{i}$ are independent of $\eps$. 
\end{assumption}

\begin{remark}[consequences of the assumptions]\label{remConseq}
\quad 
\begin{enumerate}[(i)]
\item Assumption~\ref{AssPattern} implies that $\F_{\rm per}=\bigcup_{k\in\Z^d}(\F+k)$ is connected, hence $A^\ast$ is symmetric positive definite. Indeed, $\xi\cdot A^\ast\xi=\frac{1}{|\F|}\int_\F|\nabla w_\xi|^2$ with $w_\xi=\xi_iw_i$, $\pi_\xi=\xi_i\pi_i$. If this vanishes, then $w_\xi=0$ and $\nabla\pi_\xi=\xi$ on the connected set $\F_{\rm per}$, so $\pi_\xi(y)=\xi\cdot y+c$, which is periodic only if $\xi=0$. Assumption~\ref{AssPattern} covers isolated Lipschitz inclusions $\bar B\subset Y$ in any dimension and, for $d=3$, connected solid structure crossing the faces of $Y$, provided the fluid remains connected. It excludes tangential contacts between solid parts, which create cusps in the fluid.
\item Assumption~\ref{AssOmeps}(c) allows for a solid attached to the wall $\partial\Omega$ but forbids tangential contacts or similar situations creating small features of length scale $<<\eps$. The configuration where a strip of width $\sim\eps$ without any solid in it, assumed in previous studies \cite{StokesPart2, shen22, balazi24}, is also allowed. Assumption~\ref{AssOmeps}(c) is essentially void in this case, except for the careful choice of the mesh cells $\Tbl$ in the boundary layer, which is always possible. 
\end{enumerate}
\end{remark}

\paragraph{Interpolation operator.} Let $\Oin=\big(\bigcup_{k\in\I_0}
\overline{Y^\eps_k}\big)^\circ$. For a node $a\in\eps\Z^d$, let $\varphi_a$ be the continuous $\mathbb Q_1$ hat function of $a$, i.e. $\varphi_a(x)=\prod_{s=1}^d\phi((x_s-a_s)/\eps)$ with $\phi:\R\to[0,1]$ the piecewise affine function on the grid $\Z$ such that $\phi(t)=\delta_{t0}$ at all nodes $t\in\Z$. Let $\omega_a$ be support of $\varphi_a$,
i.e. the union of the $2^d$ cells of the grid $\{Y^\eps_k\}_{k\in\Z^d}$ having
$a$ as a vertex. If $a\in\Oin$, then $\omega_a\subset\Oin$. For $v\in L^1(\Omega)$ we define the Clément-type interpolant 
\begin{equation}\label{defI}
	I_\eps^0v=\sum_{a\in\eps\Z^d\cap\Oin}\Big(\frac{1}{|\omega_a|}
	\int_{\omega_a}v\Big)\,\varphi_a .
\end{equation}
$I_\eps^0v$ is continuous on $\bar{\Omega}$, a $\mathbb Q_1$ polynomial on every grid cell in $\Tbl$, and vanishes
outside $\Oin$, so
$I_\eps^0v\in W^{1,\infty}(\R^d)$. 

\begin{theorem}\label{thmMain}
Let Assumptions~\ref{AssPattern} and \ref{AssOmeps} hold and $\eps\leqslant 1$. Let $(u^\eps,p^\eps)\in H^1_0(\Omega^\eps)^d\times L^2_0(\Omega^\eps)$ solve \eqref{genP}--\eqref{genB}, and let $p^\ast\in H^1(\Omega)$ be the weak solution of \eqref{darcy}--\eqref{DarcyB}, normalized by $\int_{\Omega^\eps}p^\ast=0$. Assume $f-\nabla p^\ast\in H^1(\Omega)^d$. Then, with $\hat u^\ast$ defined in \eqref{reghomvel},
\begin{equation}\label{mainest}
  \frac{1}{\eps^2}\|u^\eps-\hat u^\ast\|_{L^2(\Omega^\eps)}+\frac1\eps|u^\eps-\hat u^\ast|_{H^1(\Omega^\eps)}+\|p^\eps-p^\ast\|_{L^2(\Omega^\eps)}\leqslant C\sqrt\eps\,\|f-\nabla p^\ast\|_{H^1(\Omega)} 
\end{equation}
with a constant $C$ independent of $\eps$. $C$ depends only on $\Omega$, $\F$, and the constants in Assumption \ref{AssOmeps}.
\end{theorem}

\begin{remark}[on the regularity hypothesis]
	If $\Omega$ is a domain of class $C^{1,1}$, and $f\in H^1(\Omega)^d$, then $p^\ast\in H^2(\Omega)$ \cite{grisvard} and $\|f-\nabla p^\ast\|_{H^1}\leqslant C\|f\|_{H^1}$. Hence, the error bound (\ref{mainest}) holds assuming only $f\in H^1(\Omega)^d$.
	On a general Lipschitz domain, $f-\nabla p^\ast\in H^1(\Omega)^d$ is an additional hypothesis on $(f,p^\ast)$: the Neumann problem \eqref{darcy}--\eqref{DarcyB} only gives $p^\ast\in H^{3/2}(\Omega)$ in general. 
\end{remark}

\begin{remark}[error without boundary layer]
If the boundary layer is absent, i.e. if $\Omega$ is a box, consisting of whole $\eps$-cells $Y^\eps_k$, with periodic conditions on $\partial\Omega$, and assuming that $\hat u^\ast$ is based on the interpolant $I_\eps$ with the sum in (\ref{defI}) extended to all the nodes $\eps\Z^d$”, the same proof as below gives 
$$
  \frac{1}{\eps^2}\|u^\eps-\hat u^\ast\|_{L^2(\Omega^\eps)}+\frac1\eps|u^\eps-\hat u^\ast|_{H^1(\Omega^\eps)}+\|p^\eps-p^\ast\|_{L^2(\Omega^\eps)}\leqslant C\eps\|f-\nabla p^\ast\|_{H^1(\Omega)} 
  \leqslant C\eps\|f\|_{H^1(\Omega)}.$$ 
It suffices to get rid of all the bounds related to the boundary layer, no longer needed, and to adapt the interpolation error estimates (\ref{I1})--(\ref{I4}) gaining there an extra $\sqrt\eps$ (this is seen immediately in the proof of Lemma \ref{lemInterp}).
\end{remark}

\section{Auxiliary results}\label{secLemmas}

This section collects several Lemmas necessary for the proof of the main theorem in the forthcoming Section \ref{secProofThM}. We shall systematically use the letter $C$ to denote constants depending only on $d$, $\F$, $\Omega$ and the constants in Assumption~\ref{AssOmeps}; they never depend on $\eps$ or on the data. We assume $\eps\leqslant1$.  

\paragraph{Compact notation.} For a $\Z^d$-periodic function $\phi$, we systematically write $(\phi)_\eps(x)=\phi(x/\eps)$. For a measurable set $D$, $\|\cdot\|_D$ denotes the norm in $L^2(D)$. Summation over repeated indices $i$ is implied.

\paragraph{Notations related to the $\eps$-mesh and fluid parts.} We recall the notion of the mesh $\T$ and its submeshes $\Tin$, $\Tbl$ from Assumption~\ref{AssOmeps}. A {face} is a set $F=T\cap T'$ with $T\ne T'\in\T$ such that the $(d-1)$-measure $|F|$ is positive ($|F|>0$). We denote by $\E$ the set of faces and by $\Ein$ the set of faces between two interior cells $T\ne T'\in\Tin$. For every face we fix once and for all an orientation, i.e. we call one of the two cells, sharing $F$, $T$ and the other $T'$, and denote by $n_F$ the unit normal to $F$ pointing from $T$ to $T'$. \\
For any cell $T\in\T$ we shall write $\F_T^\eps$ for its fluid part, i.e. $\F_T^\eps=T\cap\Omega^\eps$ (in particular, for the interior cells $T=\overline{Y^\eps_k}$, $k\in\I_0$, the fluid part is $\F_T^\eps=\F_k^\eps$). Similarly, for a face $F=T\cap T'$ between two cells $T, T' \in \T$, we shall denote by $\F_F^\eps$ the fluid part of the union of these cells, i.e. $\F_F^\eps=(T\cup T')^\circ\cap\Omega^\eps$.

\subsection{Poincaré and trace inequalities}

\begin{lemma}[strip estimates on a Lipschitz domain]\label{lemStrip}
	For $t>0$, let $\Omega_t=\{x\in\Omega:\dist(x,\partial\Omega)<t\}$ be the boundary strip of width $t$. 
\begin{enumerate}[(a)]
\item if $w\in H^1(\Omega)$, then $\|w\|_{\Omega_t}\leqslant C\sqrt{t}\|w\|_{H^1(\Omega)}$;
\item if $w\in H^1_0(\Omega)$, then $\|w\|_{\Omega_t}\leqslant Ct\|\nabla w\|_{\Omega_{Ct}}$.
\end{enumerate}
\end{lemma}

\begin{proof}
	Part (a) is borrowed from \cite[Proposition 3.1.7]{shen2018}. There, the statement is in fact stronger: $\|w\|_{\Omega_t}\leqslant C\sqrt{t}\|w\|_{W^{1,q}(\Omega)}$ with $q=\frac{2d}{d+1}$. Our statement follows from this since $q<2$. 
	
To prove part (b), we first consider a function $v\in H^1(Q^+_r)$, where
$Q^+_r=\{x=(x',x_d)\in\R^d:|x'|<r,\ x_d>0\}$, vanishing on $\{x_d=0\}$. For
$s>0$, 
$
\int_{|x'|<r}|v(x',s)|^2\,dx'\leqslant2\int_0^s\int_{|x'|<r}|v|\,|\partial_dv|\,dx'\,dx_d .
$
Integrating over $s\in(0,t)$ gives
\[
\int_0^t\int_{|x'|<r}|v|^2\,dx'\,dx_d\leqslant2t\int_0^t\int_{|x'|<r}|v|\,|\nabla v|\,dx'\,dx_d ,
\]
and, by the Cauchy--Schwarz inequality,
\[
\Big(\int_0^t\int_{|x'|<r}|v|^2\,dx'\,dx_d\Big)^{1/2}\leqslant2t\Big(\int_0^t\int_{|x'|<r}|\nabla v|^2\,dx'\,dx_d\Big)^{1/2}.
\]
Part (b) follows by covering $\partial\Omega$ with coordinate patches in which $\Omega$ lies above the graph of a Lipschitz function $\gamma$, and applying the result above after the bi-Lipschitz change of variables
$(x',x_d)\mapsto(x',x_d-\gamma(x'))$, which maps $\Omega_t$ into a strip of width $Ct$.
\end{proof}

\begin{lemma}[Poincar\'e inequality]\label{lemPoin}
	For $\phi\in H^1_0(\Omega^\eps)$, $\|\phi\|_{\Omega^\eps}\leqslant C\eps|\phi|_{H^1(\Omega^\eps)}$.
\end{lemma}
\begin{proof}
	Extend $\phi$ by $0$ to $\Omega$. By Assumption~\ref{AssOmeps}, $\Omega\subset\overline{\Omega_{c_0\eps}}\cup\bigcup_{k\in\I_0}\overline{Y^\eps_k}$ in the notations of Lemma~\ref{lemStrip}. We can apply part (b) of this Lemma to bound $\|\phi\|_{\Omega_{c_0\eps}}\leqslant C\eps|\phi|_{H^1(\Omega_{C\eps})}$. On $Y_k^\eps$, $k\in\I_0$, $\phi$ vanishes on $B_k^\eps$, and the scaled Poincar\'e inequality on $Y$ for functions vanishing on the non-empty open set $B$ gives $\|\phi\|_{Y_k^\eps}\leqslant C\eps|\phi|_{H^1(Y_k^\eps)}$.
\end{proof}

\begin{lemma}[trace inequality on the cells]\label{lemTrace}
For every $T\in\T$ and $w\in H^1(T^\circ)$,
\[
  \|w\|^2_{\partial T}\leqslant C\big(\eps^{-1}\|w\|^2_T+\eps|w|^2_{H^1(T)}\big).
\]
Consequently, for $w\in H^1(\Omega)$, $\sum_{F\in\E}\|w\|^2_F\leqslant C\big(\eps^{-1}\|w\|^2_\Omega+\eps|w|^2_{H^1(\Omega)}\big)$.
\end{lemma}

\begin{proof}
Lemma~\ref{lemTrace} follows by scaling from the trace inequality
$\|w\|^2_{\partial\hat T}\leqslant C(\|w\|^2_{\hat T}+\|w\|_{\hat T}\|\nabla w\|_{\hat T})$
on the Lipschitz domain $\hat T=\eps^{-1}T$, whose constant depends only on the
Lipschitz character of $\hat T$: it is obtained by integrating
$\Div(\mu w^2)$ with a smooth vector field $\mu$ transversal to
$\partial\hat T$, constructed from the local graph representations of
$\partial\hat T$, see \cite[Lemma~1.5.1.9 and Theorem~1.5.1.10]{grisvard}.
\end{proof}

\subsection{Pressure means and the inf-sup condition}

\begin{lemma}[Bogovski\u{\i} operator]\label{lemBog}
Let $\omega\subset\R^d$ be a John domain. For every $q\in L^2(\omega)$ with $\int_\omega q=0$ there is $v\in H^1_0(\omega)^d$ with $\Div v=q$ and $|v|_{H^1(\omega)}\leqslant C_B\|q\|_\omega$, where the constant $C_B$ depends only on the John constant $J$ of $\omega$ as in Definition \ref{deftnqohn}.
\end{lemma}

\begin{proof}
The operator defining $v$ is constructed in \cite{bogovskii}. The fact that this construction is valid on John domains and that $C_B$ depends only on $J$ is proved in \cite{acosta06}. 
\end{proof}

For $q\in L^2(\Omega^\eps)$ and $T\in\T$ let $q_T=|\F_T^\eps|^{-1}\int_{\F_T^\eps}q$ and $\Pi_\eps q$ be the piecewise constant function on the mesh $\T$ defined on every $T\in\T$ by $\Pi_\eps q|_T=q_T$. $\Pi_\eps$ maps $L^2_0(\Omega^\eps)$ into itself (after restricting $\Pi_\eps q$ to the fluid part $\F_T^\eps$ on every cell $T$). We introduce the jumps $\jump{\Pi_\eps q}_F=q_T-q_{T'}$ across any face $F\in\E$,  $F=T\cap T'$, using the conventions adopted in the beginning of this Section. Set
\begin{equation}\label{deftnq}
  \tn{q}^2=\eps\sum_{F\in\E}\|\jump{\Pi_\eps q}\|^2_F+\|(I-\Pi_\eps)q\|^2_{\Omega^\eps} .
\end{equation}

\begin{lemma}[discrete Poincar\'e inequality]\label{lemDiscPoin}
For every $q\in L^2_0(\Omega^\eps)$,
\begin{equation}\label{discPoin}
  \|\Pi_\eps q\|_{\Omega^\eps}\leqslant C\eps^{-1}\tn{q}\qquad\text{and}\qquad\|q\|_{\Omega^\eps}\leqslant C\eps^{-1}\tn{q}.
\end{equation}
\end{lemma}

\begin{proof}
	Let $Q\in L^2(\Omega)$ be equal to $\Pi_\eps q$ on $\Omega^\eps$ and to $0$ in
	the solid. Then $\int_\Omega Q=\sum_{T\in\T}|\F_T^\eps|q_T=\int_{\Omega^\eps}q=0$ and
	$\|Q\|_\Omega=\|\Pi_\eps q\|_{\Omega^\eps}$. Let $V\in H^1_0(\Omega)^d$ satisfy $\Div V=Q$ and $\|V\|_\Omega+|V|_{H^1(\Omega)}\leqslant
	C\|Q\|_\Omega$ (this $V$ is given by Lemma \ref{lemBog} since the Lipschitz domain $\Omega$ is also John, and the constant $C$ depends only on $\Omega$ through its John constant and its diameter, involved in the Poincaré-Friedrichs inequality). Observe
	$$		\|\Pi_\eps q\|^2_{\Omega^\eps}
	=\sum_{T\in\T}q_T\int_{\F_T^\eps}q_T
	=\sum_{T\in\T}q_T\int_T\Div V
	=\sum_{F\in\E}\jump{\Pi_\eps q}_F\,V\cdot n_F\,.
	$$
	The last equality follows from the divergence theorem on every cell $T$, rearranging the sums over the cell boundaries, and the fact that $V\cdot n$ vanishes on $\partial\Omega$, so that only the contributions on faces $F\in\E$ survive. We now use	 
	the Cauchy--Schwarz inequality and Lemma~\ref{lemTrace}:
	\begin{multline*}
		\|\Pi_\eps q\|^2_{\Omega^\eps} \leqslant\Big(\eps\sum_{F\in\E}\|\jump{\Pi_\eps q}\|^2_F\Big)^{1/2}
		\Big(\eps^{-1}\sum_{F\in\E}\|V\|^2_F\Big)^{1/2}\\
		\leqslant C\tn{q}\big(\eps^{-2}\|V\|_\Omega^2+|V|_{H^1(\Omega)}^2\big)^{1/2}
		\leqslant C\eps^{-1}\tn{q}\,\|\Pi_\eps q\|_{\Omega^\eps}
	\end{multline*}
	since $\eps<1$. This proves the first inequality in \eqref{discPoin}. The second inequality follows from
	$\|q\|_{\Omega^\eps}^2=\|\Pi_\eps q\|_{\Omega^\eps}^2+\|(I-\Pi_\eps)q\|_{\Omega^\eps}^2$ and $\eps<1$.
\end{proof}

\begin{lemma}[velocity-pressure inf-sup]\label{lemInfSup}
For every $p\in L^2_0(\Omega^\eps)$ there is $w\in H^1_0(\Omega^\eps)^d$ with
\[
  (p,\Div w)_{\Omega^\eps}=\tn{p}^2,\qquad\beta|w|_{H^1(\Omega^\eps)}\leqslant\tn{p},
\] with $\beta\in(0,1]$ independent of $\eps$.
\end{lemma}

\begin{proof}
For $F\in\E$ with cells $T,T'$, Lemma~\ref{lemBog} gives $\Phi_F\in H^1_0(\F_F^\eps)^d$ with
\begin{equation}\label{PhiF}
  \Div\Phi_F=\frac{\mathbf 1_{\F_T^\eps}}{|\F_T^\eps|}-\frac{\mathbf 1_{\F_{T'}^\eps}}{|\F_{T'}^\eps|},\qquad|\Phi_F|_{H^1(\F_F^\eps)}\leqslant C_B\big(|\F_T^\eps|^{-1}+|\F_{T'}^\eps|^{-1}\big)^{1/2}\leqslant C\eps^{-d/2},
\end{equation}
using $|\F_T^\eps|\geqslant c\eps^d$ (Assumption~\ref{AssOmeps}(c)). Domain $\F_F^\eps$ is indeed John with constant $J$ by Assumption~\ref{AssOmeps}(c), if one of the cells $T,T'$ is in $\Tbl$. For patches $\F_F^\eps$ composed of interior cells, the John constant is also independent of $\eps$, cf. Assumption~\ref{AssPattern}, since the John constant is scale invariant. 

Set $v_p=\eps\sum_{F\in\E}|F|(p_T-p_{T'})\Phi_F$. Since $(p,\Div\Phi_F)_{\Omega^\eps}=p_T-p_{T'}$, we get $$(p,\Div v_p)_{\Omega^\eps}=\eps\sum_{F\in\E}|F|(p_T-p_{T'})^2=\eps\sum_{F\in\E}\|\jump{\Pi_\eps p}\|^2_F\,.$$ The supports of the $\Phi_F$ overlap boundedly and $|F|\leqslant C\eps^{d-1}$, so
\begin{equation}\label{vpbound}
  |v_p|^2_{H^1(\Omega^\eps)}\leqslant C\eps^2\sum_{F\in\E}|F|^2(p_T-p_{T'})^2\eps^{-d}\leqslant C\eps\sum_{F\in\E}\|\jump{\Pi_\eps p}\|^2_F .
\end{equation}
%Moreover, $\Div v_p$ is constant on every $\F_T^\eps$, hence orthogonal to $(I-\Pi_\eps)p$.

For every $T\in\T$, $(I-\Pi_\eps)p$ is of zero mean on $\F_T^\eps$. Thus, Lemma~\ref{lemBog} gives $v_T\in H^1_0(\F_T^\eps)^d$ with $\Div v_T=(I-\Pi_\eps)p$ on $\F_T^\eps$ and $|v_T|_{H^1(\F_T^\eps)}\leqslant C_B\|(I-\Pi_\eps)p\|_{\F_T^\eps}$. The constant $C_B$ is $\eps$-independent by the same considerations as above. 

Let $v_{p0}=\sum_{T\in\T}v_T$ (each extended by zero) and $w=v_p+v_{p0}$. Then $$(p,\Div w)_{\Omega^\eps}=(p,\Div v_{p})_{\Omega^\eps}+(p-\Pi_\eps p,\Div v_{p0})_{\Omega^\eps}=(p,\Div v_{p})_{\Omega^\eps}+\|(I-\Pi_\eps)p\|_{\Omega^\eps}^2=\tn{p}^2,$$ since $\Div v_{p0}$ has zero mean on every $\F_T^\eps$. Recalling  $|v_{p0}|_{H^1(\Omega^\eps)}\leqslant C_B\|(I-\Pi_\eps)p\|_{\Omega^\eps}$ and combining with (\ref{vpbound}) gives the announced bound on the norm of $w$. 
\end{proof}

\begin{remark}\label{remtnqOld}
	An alternative definition of the mesh-dependent norm on the pressure space could be
\begin{equation}\label{deftnqOld}  
	\interleave q \interleave_{\eps,\textrm{alt}}^2
	=\eps^2\sum_{T\in\T}\|\Pi_\eps q\|^2_T+\|(I-\Pi_\eps)q\|^2_{\Omega^\eps} .
\end{equation}
With this definition, Lemma \ref{lemDiscPoin} would be a triviality, and Lemma \ref{lemInfSup} would be simpler to prove. One can work with such definition as well, and obtain almost the same results as with the present definition (\ref{deftnq}), but the regularity requirement would be stronger (typically $f-\nabla p^\ast \in H^{3/2}$) and the treatment of the boundary layer contributions would be more involved.
\end{remark}

We shall also use the combined velocity-pressure space $\X=H^1_0(\Omega^\eps)^d\times L^2_0(\Omega^\eps)$, equipped with the norm $\tn{(v,q)}$ defined for any $(v,q)\in\X$ by
\begin{equation}\label{defnorm}
	\tn{(v,q)}^2=|v|^2_{H^1(\Omega^\eps)}+\tn{q}^2 .
\end{equation}
Introduce the combined bilinear form $c:\X\times\X\to\R$ associated to the Stokes equations by
 $$c((u,p),(v,q))=(\nabla u,\nabla v)_{\Omega^\eps}-(p,\Div v)_{\Omega^\eps}-(q,\Div u)_{\Omega^\eps}.$$
\begin{lemma}[inf-sup for the combined bilinear form $c$]\label{lemInfSupc}
For every $(u,p)\in\X$, there exists $(v,q)\in\X\setminus\{0\}$ such that
\begin{equation}\label{infsup}
	\frac{c((u,p),(v,q))}{\tn{(v,q)}}\geqslant\beta_\ast\tn{(u,p)}
\end{equation}
with $\beta_\ast=\beta^2/(2\sqrt3)$ and $\beta$ from  Lemma \ref{lemInfSup}.
\end{lemma}
\begin{proof}
Take any $(u,p)\in\X$  and choose $w$ as in Lemma \ref{lemInfSup}. Then
\[
c((u,p),(-w,0))=-(\nabla u,\nabla w)_{\Omega^\eps}+\tn{p}^2\geqslant-\frac{1}{\beta}|u|_{H^1(\Omega^\eps)}\tn{p}+\tn{p}^2,\qquad c((u,p),(u,-p))=|u|_{H^1(\Omega^\eps)}^2 .
\]
With $v=u-\beta^2w$ and $q=-p$, we get $$c((u,p),(v,q))\geqslant |u|_{H^1(\Omega^\eps)}^2-\beta|u|_{H^1(\Omega^\eps)}\tn{p}+\beta^2\tn{p}^2\geqslant\frac12|u|_{H^1(\Omega^\eps)}^2+\frac{\beta^2}2\tn{p}^2\geqslant\frac{\beta^2}{2}\tn{(u,p)}^2$$ and $\tn{(v,q)}^2 \leqslant (|u|_{H^1(\Omega^\eps)}+\beta\tn{p})^2+\tn{p}^2 \leqslant 3\tn{(u,p)}^2$.
\end{proof}

\subsection{Interpolation and the flux corrector}

\begin{lemma}\label{lemInterp}
For $v\in H^1(\Omega)$, with $I_\eps^0$ defined in \eqref{defI},
\begin{align}
  \|I_\eps^0v\|_{\Omega}&\leqslant C\|v\|_{\Omega},\label{I0}\\
  \|v-I_\eps^0v\|_{\Omega}&\leqslant C\sqrt\eps\,\|v\|_{H^1(\Omega)},\label{I1}\\
  |I_\eps^0v|_{H^1(\Omega)}&\leqslant\frac{C}{\sqrt\eps}\|v\|_{H^1(\Omega)},\label{I2}\\
  \Big(\sum_{F\in\Ein}\|\nabla I_\eps^0v\|^2_F\Big)^{1/2}&\leqslant\frac{C}{\eps}\|v\|_{H^1(\Omega)},\label{I3}\\
  \Big(\sum_{F\in\E}\|v-I_\eps^0v\|^2_F\Big)^{1/2}&\leqslant C\|v\|_{H^1(\Omega)} .\label{I4}
\end{align}
\end{lemma}

\begin{proof}
By scaling and \eqref{defI}, we have for any $T \in \Tin$,
\begin{equation}
	\label{stabT} \|I^0_{\eps} v\|_T^2 \leqslant C \eps^d \sum_a | I^0_{\eps} v
	(a) |^2 \leqslant C \sum_a \|v\|^2_{\omega_a} \leqslant C
	\|v\|^2_{\omega_T}, \qquad |I^0_{\eps} v|_{H^1 (T)} \leqslant \frac{C}{\eps}
	\|I^0_{\eps} v\|_T \leqslant \frac{C}{\eps} \|v\|_{\omega_T},
\end{equation}
where the sum is taken over all the vertices of $T$ inside $\Oin$, and
$\omega_T$ denotes the union of node-associated patches $\omega_a$ for all
such vertices $a$. Bound \eqref{I0} follows by summing the first inequality in
(\ref{stabT}) over all $T \in \Tin$, since the patches $\omega_T$ overlap
boundedly.

Let $\mathcal{T}^{\eps}_{\mathrm{in}, 0}$ be the set of cells $T \in \Tin$
whose $2^d$ vertices all lie inside $\Oin$, and
$\mathcal{T}^{\eps}_\textrm{in,bl}= \Tin \setminus
\mathcal{T}^{\eps}_{\mathrm{in}, 0}$. If $T \in
\mathcal{T}^{\eps}_{\mathrm{in}, 0}$, then $\omega_T$ is a cube of side $3
\eps$ and $I^0_{\eps} c = c$ on $T$ for any constant $c$. Taking $c$ as the
mean of $v$ on $\omega_T$, \eqref{stabT} and the Poincar{\'e} inequality on
$\omega_T$ give, on $T \in \mathcal{T}^{\eps}_{\mathrm{in}, 0}$,
\begin{multline}\label{InStabBetter} 
	\|v - I^0_{\eps} v\|_T = \|v - c - I^0_{\eps}  (v - c)
	\|_T \leqslant C \|v - c\|_{\omega_T} \leqslant C \eps \| \nabla
	v\|_{\omega_T}, \\ |I^0_{\eps} v|_{H^1 (T)} = |I^0_{\eps} (v - c) |_{H^1
		(T)} \leqslant C \| \nabla v\|_{\omega_T} .
\end{multline}
If $T \in \mathcal{T}^{\eps}_\textrm{in, bl}$, then $T$
touches $\partial \Oin$, which lies at distance at most $c_0 \eps$ from
$\partial \Omega$ by Assumption~\ref{AssOmeps}(b). Hence $T$ and $\omega_T$,
as well as $\Obl = \Omega \setminus \Oin$, are contained in the strip
$\Omega_{C \eps}$. Summing over all cells and recalling $I^0_{\eps} v = 0$ on
$\Obl$, we obtain thanks to Lemma~\ref{lemStrip}(a), and
(\ref{stabT})--(\ref{InStabBetter}),
\begin{equation*}
	\|v - I^0_{\eps} v\|^2_{\Omega}  = \sum_{T \in
		\mathcal{T}^{\eps}_{\mathrm{in}, 0}} \|v - I^0_{\eps} v\|^2_T + \sum_{T \in
		\mathcal{T}^{\eps}_\mathrm{in,bl}} \|v - I^0_{\eps} v\|^2_T +
	\|v\|^2_{\Obl} \leqslant C \left( \eps^2 |v|^2_{H^1 (\Omega)}
	+\|v\|^2_{\Omega_{C \eps}} \right) \leqslant C \eps \|v\|^2_{H^1 (\Omega)},
\end{equation*}
and similarly
\begin{equation*}
	|I^0_{\eps} v|^2_{H^1 (\Omega)}  \leqslant C \left( |v|^2_{H^1 (\Omega)} +
	\eps^{- 2} \|v\|^2_{\Omega_{C \eps}} \right) \leqslant \frac{C}{\eps}
	\|v\|^2_{H^1 (\Omega)},
\end{equation*}
which are \eqref{I1} and \eqref{I2}.
	
	It remains to prove the estimates on faces. For \eqref{I3}, a scaling argument on	a cube $Y_k^{\eps}$ gives $\| \nabla I^0_{\eps} v\|^2_F \leqslant C \eps^{- 1} |I^0_{\eps}
	v|^2_{H^1 (Y_k^{\eps})}$ for a face $F$ of $Y_k^{\eps}$, $k\in\I_0$. Thus, the left-hand side of (\ref{I3}) is bounded by $C \eps^{- 1/2} |I^0_{\eps}v|_{H^1 (\Oin)}$, and (\ref{I3}) follows from 	\eqref{I2}. For \eqref{I4}, Lemma~\ref{lemTrace} gives $\sum_{F \in \E} \|v
	- I^0_{\eps} v\|^2_F \leqslant C (\eps^{- 1} \|v - I^0_{\eps} v\|^2_{\Omega} +
	\eps |v - I^0_{\eps} v|^2_{H^1 (\Omega)})$, so that \eqref{I4} follows from \eqref{I1} and	\eqref{I2}.
\end{proof}

We now recall the classical flux corrector of periodic homogenization \cite{jko,shen2018}.

\begin{lemma}[flux corrector]\label{lemFlux}
For $i=1,\dots,d$ there is a $\Z^d$-periodic skew-symmetric matrix field $\phi^{(i)}\in H^2_{\rm per}(\R^d)^{d\times d}$ such that
$${w}_i - | \F | A^{\ast} e_i = \Div  \phi^{(i)}, \quad \text{with} \quad
\left( \Div  \phi^{(i)} \right)_a = \sum_{b = 1}^d \partial_b
\phi^{(i)}_{ba}$$
We assume here that $w_i$ is extended by 0 in $B=Y\setminus\F$.
\end{lemma}

\begin{proof}
Since $w_i \in H^1(\F)$, $\Div w_i = 0$ in $\F_{\mathrm{per}}$ and $w_i = 0$ on
$\partial \F_{\mathrm{per}}$, extending ${w}_i$ by 0 in $B$ gives a field in $H^1_{\mathrm{per}} (\R^d)^d$ and its divergence vanishes everywhere. Its mean over
$Y$ is $\int_{\F} w_i = | \F | A^{\ast} e_i$. Let $\psi \in H^3_{\mathrm{per}}
(\R^d)^d$ solve $- \Delta \psi = {w}_i - | \F | A^{\ast} e_i$ and set
$\phi_{ba}^{(i)} = \partial_a \psi_b - \partial_b \psi_a$. Then $(\Div
\phi^{(i)})_a = \partial_a \Div \psi - \Delta \psi_a$. Observe $- \Delta
\Div \psi = \Div {w}_i = 0$. Thus, $\Div \psi$ is a periodic
harmonic function, hence constant, and $\Div \phi^{(i)} = - \Delta \psi $.\end{proof}

\begin{corollary}\label{corfluxid}
	Let $h$ be a function of class $C^2$ on a grid cell  $Y^\eps_k$. Then 
	\begin{equation}\label{fluxid}
		\eps^2\int_{Y^\eps_k}( w_i)_\eps\cdot\nabla h=\eps^2|\F|\int_{Y^\eps_k}A^\ast e_i\cdot\nabla h+\eps^3\int_{\partial Y^\eps_k}(\phi^{(i)})_\eps\,\nabla_\tau h\cdot n ,
	\end{equation}
	where $\nabla_\tau$ denote the tangential part of the gradient on $\partial Y^\eps_k$.
\end{corollary}
\begin{proof}
We write  $  ( w_i)_\eps=|\F|A^\ast e_i + \eps\Div\big(\phi^{(i)}\big)_\eps$, and perform the integration by parts in the integral
\[\eps^3\int_{Y^\eps_k}\Div\big(\phi^{(i)}\big)_\eps\cdot\nabla h
=\eps^3\int_{\partial Y^\eps_k}(\phi^{(i)})_\eps\,\nabla h\cdot n 
-\eps^3\int_{Y^\eps_k}\sum_{a,b=1}^d(\phi^{(i)}_{ab})_\eps\,\partial_a\partial_bh \,.
\] 
The last integral vanishes since $\phi$ is a skew-symmetric matrix and the matrix of second derivatives of $h$ is symmetric. Writing, on $\partial Y^\eps_k$, $\nabla h = (\nabla h \cdot n)n + \nabla_\tau h$, we conclude $\nabla h\cdot(\phi^{(i)})_\eps\,n=(\phi^{(i)})_\eps\,\nabla_\tau h\cdot n$ since $(\phi^{(i)})_\eps\,n\cdot n=0$ by skew-symmetry of $\phi^{(i)}$.
\end{proof}

\begin{remark}[products with oscillating functions] \label{remprod} 	
	We shall frequently encounter the products of $\eps$-periodic functions, eg. $w_i,\nabla w_i,\pi_i$ rescaled by $\eps$, with $\mathbb Q_1$ functions on mesh $\Tin$. If one needs to calculate the norm of such a product, then multiplication by the periodic function will act as multiplication by a constant. Indeed, let $\psi\in L^2(\F)$ be $\Z^d$-periodic  and $h$ be a polynomial of degree up to $d$ on $Y_k^\eps$. Scaling estimates give then
	\begin{equation*}
		\|(\psi)_\eps h\|_{\F^\eps_k}\leqslant\eps^{d/2}\|\psi\|_{\F}\|h\|_{L^\infty(Y^\eps_k)}\leqslant C\|\psi\|_{\F}\|h\|_{Y^\eps_k}\leqslant C\|h\|_{Y^\eps_k}.
	\end{equation*}
	In the last inequality we made the generic constant $C$ absorb the norm $\|\psi\|_{\F}$ since, in all the cases encountered, it will depend only on $\F$, and the constants $C$ are allowed to hide this dependence.   
	%The same observation is valid for the $L^2$ norms on faces in $\Ein$. This will be useful for the products with $\phi^{(i)}$.
\end{remark}

\section{Proof of Theorem~\ref{thmMain}}\label{secProofThM}

Let us denote $\bar{u} = f - \nabla p^{\ast}$ and recall Definition (\ref{eq:homvel}) of $\hat{u}^{\ast}$, which we rewrite as $\hat{u}^{\ast} =
\eps^2 ({w}_i)_{\eps}  \, I_{\eps}^0 \bar{u}$ and introduce
the pressure reconstruction
\begin{equation}
	\label{defansatz} p^{\eps, 1} = p^{\ast} + \eps ({\pi}_i)_{\eps} 
	\, I_{\eps}^0  \bar{u}_i - m_{\eps}, \qquad m_{\eps} =
	\frac{1}{| \Omega^{\eps} |}  \int_{\Omega^{\eps}} \eps
	({\pi}_i)_{\eps}  \, I_{\eps}^0  \bar{u}_i .
\end{equation}
Since $I_{\eps}^0  \bar{u}$ vanishes outside $\Oin$, $\hat{u}^{\ast} \in H^1_0
(\Omega^{\eps})^d$. By construction, $p^{\eps, 1} \in L^2_0
(\Omega^{\varepsilon})$. Take some $z \in H^1_0 (\Omega^{\eps})^d$, to be
chosen later in the proof. \ We have $(u^{\eps} - \hat{u}^{\ast} + z, p^{\eps}
- p^{\eps, 1}) \in \X$. By Lemma~\ref{lemInfSupc} there is $(v, q) \in \X$
with $\tn{(v, q)} = 1$ and
\begin{multline}\label{erreq} 
	\beta_{\ast} \tn{(u^{\eps} - \hat{u}^{\ast} + z, p^{\eps} -
		p^{\eps, 1})} \leqslant c ((u^{\eps} - \hat{u}^{\ast} + z, p^{\eps} - p^{\eps,	1}), (v, q)) \\
		= \underbrace{(\nabla (u^{\eps} - \hat{u}^{\ast} + z), \nabla
		v)_{\Omega^{\eps}}}_I - \underbrace{(p^{\eps} - p^{\eps, 1}, \Div
		v)_{\Omega^{\eps}}}_{II} - \underbrace{(q, \Div (u^{\eps} - \hat{u}^{\ast} +
		z))_{\Omega^{\eps}}}_{III},
\end{multline}
We shall now treat the terms on the right-hand side of (\ref{erreq}) and
choose $z$, finally obtaining the announced error estimate.

\paragraph{The terms $I - II$.} The weak form of \eqref{genP} and $f = \bar{u}
+ \nabla p^{\ast}$ give
\[ (\nabla u^{\eps}, \nabla v)_{\Omega^{\eps}} = (\bar{u}, v)_{\Omega^{\eps}} - (p^{\ast}, \Div v)_{\Omega^{\eps}} +
(p^{\eps}, \Div v)_{\Omega^{\eps}} . \]
We are now going to write out the analogue of this relation for
$\hat{u}^{\ast}$:\footnote{We treat the gradients of vector fields as matrices following the convention $(\nabla g)_{ab}=\partial_ag_b$ for a vector field $g=(g_b)_{b=1,\ldots,d}$. All the vectors are assumed to be columns.}
\begin{multline*}
	 (\nabla \hat{u}^{\ast}, \nabla v)_{\Omega^{\eps}} = \eps ((\nabla w_i)_{\eps}, (I^0_{\eps} 
\bar{u}_i) \nabla v)_{\Omega^{\eps}} + \eps^2 (\nabla I^0_{\eps}  \bar{u}_i, \nabla v
(w_i)_{\eps})_{\Omega^{\eps}} \\
 = \eps ((\nabla w_i)_{\eps}, \nabla (vI^0_{\eps}  \bar{u}_i))_{\Omega^{\eps}} - \varepsilon
(\nabla I^0_{\eps}  \bar{u}_i, (\nabla w_i)_{\eps} v)_{\Omega^{\eps}} + \eps^2 (\nabla
I^0_{\eps}  \bar{u}_i, (\nabla v) (w_i)_{\eps})_{\Omega^{\eps}} 
\end{multline*}
The cell problems \eqref{eq:cell}, rescaled by $\varepsilon$, read in the weak
form tested by $vI^0_{\eps}  \bar{u}_i$,
\[ \eps \left( (\nabla w_i)_{\eps}, \nabla (vI^0_{\eps}  \bar{u}_i)
\right)_{\Omega^{\eps}} - \eps \left( (\pi_i)_{\eps}, \Div (vI^0_{\eps} 
\bar{u}_i) \right)_{\Omega^{\eps}} = (I^0_{\eps}  \bar{u},
v)_{\Omega^{\eps}} \]
The test function $(I^0_{\eps}  \bar{u}_i) v$ is admissible since $I^0_{\eps} 
\bar{u}_i \in W^{1, \infty}$ and its support is in $\overline{\Oin}$, where
the perforations are periodic. Denoting $\sigma_i = \nabla w_i - \pi_i
\mathrm{I}$, we obtain
\[ (\nabla \hat{u}^{\ast}, \nabla v)_{\Omega^{\eps}} = (I^0_{\eps}  \bar{u}, v)_{\Omega^{\eps}} - \eps \left(
\nabla I^0_{\eps}  \bar{u}_i, (\sigma_i)_{\eps} v \right)_{\Omega^{\eps}} + \eps^2 (\nabla
I^0_{\eps}  \bar{u}_i, (\nabla v) (w_i)_{\eps})_{\Omega^{\eps}} + \eps \left(
(\pi_i)_{\eps} I^0_{\eps}  \bar{u}_i, \Div v \right)_{\Omega^{\eps}} . \]
Thus,
\[ I - II = (\bar{u} - I_{\varepsilon}^0  \bar{u}, v)_{\Omega^{\varepsilon}} +
\eps \left( \nabla I_{\varepsilon}^0  \bar{u}_i, (\sigma_i)_{\eps} v
\right)_{\Omega^{\varepsilon}} - \eps^2 (\nabla I^0_{\eps}  \bar{u}_i,
(\nabla v) (w_i)_{\eps})_{\Omega^{\varepsilon}} + (\nabla z, \nabla
v)_{\Omega^{\eps}} \hspace{0.17em} . \]
The constant $m_{\eps}$ in $p^{\eps, 1}$ does not contribute above, since
$\int_{\Omega^{\eps}} \Div v = 0$.          
By Lemma~\ref{lemPoin}, \eqref{I1}, \eqref{I2}, Remark \ref{remprod}, and $|v|_{H^1(\Omega^{\eps})} \leqslant
\tn{(v, q)} = 1$,
\begin{equation}
	\label{boundI} |I - II| \leqslant C \left( \eps \| \bar{u} - I^0_{\eps}  \bar{u}
	\|_{\Omega} + \eps^2 |I^0_{\eps}  \bar{u} |_{H^1 (\Omega)} \right) +
	|z|_{H^1 (\Omega^{\eps})} \leqslant C \eps^{3 / 2} \| \bar{u} \|_{H^1 (\Omega)} +
	|z|_{H^1 (\Omega^{\eps})} .
\end{equation}
\paragraph{The term $III$.} Since $\Div u^{\eps} = 0$,
\begin{equation}
	\label{term3} - III = \underbrace{(\Pi_{\eps} q, \Div
		\hat{u}^{\ast})_{\Omega^{\eps}}}_{III_1} + \underbrace{((I - \Pi_{\eps}) q, \Div
		\hat{u}^{\ast})_{\Omega^{\eps}}}_{III_2} - (q, \Div z)_{\Omega^{\eps}} .
\end{equation}
By Remark \ref{remprod}, \eqref{I2}, and $\|(I - \Pi_{\eps}) q\|_{\Omega^{\eps}} \leqslant \tn{q} = 1$,
\begin{equation}
	\label{boundIII2} |III_2 | \leqslant \eps^2 \|(I - \Pi_{\eps}) q\|_{\Omega^{\eps}} \,
	\|({w}_i)_{\eps} \nabla I^0_{\eps}  \bar{u}_i \|_{\Omega^{\eps}} \leqslant C \eps^2
	|I^0_{\eps}  \bar{u} |_{H^1 (\Omega)} \leqslant C \eps^{3 / 2} \| \bar{u} \|_{H^1
		(\Omega)} .
\end{equation}
Since $\Div \hat{u}^{\ast} = \eps^2 ({w}_i)_{\eps} \cdot \nabla
I^0_{\eps}  \bar{u}_i$ vanishes in the solid and outside $\Oin$, the identity
\eqref{fluxid} with $h = I^0_{\eps}  \bar{u}_i$ on every interior cell gives
\begin{equation}
	\label{III1p12} III_1 = \sum_{T \in \Tin} q_T \, \eps^2 \int_T
	({w}_i)_{\eps} \cdot \nabla I^0_{\eps}  \bar{u}_i = \underbrace{\eps^2
		| \F |  \sum_{T \in \T} q_T  \int_T \Div (A^{\ast} I^0_{\eps} 
		\bar{u})}_{III_1'} + \underbrace{\sum_{T \in \Tin} q_T \,
		\eps^3 \int_{\partial T} (\phi^{(i)})_{\eps} \nabla_{\tau} I^0_{\eps} 
		\bar{u}_i \cdot \, n}_{III_1''},
\end{equation}
where the sum in $III_1'$ could be extended to all cells since $I^0_{\eps} 
\bar{u} = 0$ on the layer cells.

\paragraph{The term $III_1'$.} Since $\Div (A^{\ast}  \bar{u}) = 0$ in
$\Omega$, term $III_1'$ can be rewritten as
\[ III_1' = \eps^2 | \F |  \sum_{T \in \T} q_T  \int_T \Div (A^{\ast}  \left(
I^0_{\eps}  \bar{u} - \bar{u} \right)) = \eps^2 | \F |  \sum_{T \in \T} q_T
\int_{\partial T \cap \Omega} n \cdot A^{\ast}  \left( I^0_{\eps}  \bar{u}
- \bar{u} \right) \]
The portions of the boundaries $\partial T$ lying on $\partial \Omega$ could be excluded from the last sum since \ $I^0_{\eps}  \bar{u} = 0$ and $n \cdot
A^{\ast}  \bar{u} = 0$ on $\partial \Omega$. We now rewrite the last sum over
the boundaries of cells $T \in \T$ as the sum over the faces $F \in
\E$, the Cauchy--Schwarz inequality and \eqref{I4} give

\begin{multline}\label{boundIII1}
	|III_1' | = \eps^2 | \F | \, \left| \sum_{F \in \E} \int_F
	\jump{\Pi_{\eps} q}_F \, n_F \cdot A^{\ast}  \left( I^0_{\eps} 
	\bar{u} - \bar{u} \right) \right| \\ \leqslant C \eps^2 \left( \eps \sum_{F \in \E} \|
	\jump{\Pi_{\eps} q} \|^2_F \right)^{1 / 2} \left( \eps^{- 1} \sum_{F \in \E}
	\| \bar{u} - I^0_{\eps}  \bar{u} \|^2_F \right)^{1 / 2} \leqslant C \eps^{3 / 2}
	\, \tn{q} \, \| \bar{u} \|_{H^1 (\Omega)} \leqslant C
	\eps^{3 / 2} \| \bar{u} \|_{H^1 (\Omega)}  
\end{multline}
since $\tn{q} \leqslant \tn{(v, q)} = 1$.

\paragraph{Choice of $z$.} We now choose $z$ so that to cancel exactly in
(\ref{term3}) the term $III_1''$\quad introduced in (\ref{III1p12}), i.e. so
that
\begin{equation}
	\label{zqident} III_1'' = (q, \Div z)_{\Omega^{\eps}} .
\end{equation}
We rewrite $III_1''$ as a sum over faces. For $F \in \Ein$, $F = T \cap T'$
with $T = \overline{Y^{\eps}_k}$, $T' = \overline{Y^{\eps}_l}$, set
\[ J_F = \eps^3 \int_F (\phi^{(i)})_{\eps} \nabla_{\tau} I^0_{\eps}  \bar{u}_i
\cdot \, n_F . \]
Note that $J_F$ is unambiguously defined since the tangential derivatives of
$I^0_{\eps}  \bar{u}_i$ on $F$ are the same from both sides $T$ and $T'$. On a
face of an interior cell which is not in $\Ein$, $I^0_{\eps}  \bar{u}_i$ and
its tangential derivatives vanish. Passing from the sum over the boundaries of
the interior cells to the sum over their faces, and recalling that $n_F$ is
the outer normal of $T$ and the inner normal of $T'$, we get
\[ III_1'' = \sum_{F \in \Ein} (q_T - q_{T'})  \, J_F . \]
On the other hand, the fields $\Phi_F$ defined in the proof of Lemma~\ref{lemInfSup} by (\ref{PhiF})
satisfy $(q, \Div \Phi_F)_{\Omega^{\eps}} = q_T - q_{T'}$. Hence the choice
\begin{equation}\label{defz} z = \sum_{F \in \Ein} J_F  \, \Phi_F \in H^1_0
(\Omega^{\eps})^d 
\end{equation}
gives \eqref{zqident}. To bound $z$, note that $\|(\phi^{(i)})_{\eps} \|_F \leqslant
C \eps^{(d - 1) / 2}$, so that $|J_F | \leqslant C \eps^{3 + (d - 1) / 2} \| \nabla_\tau I^0_{\eps}  \bar{u} \|_F$ by the Cauchy--Schwarz inequality. Since $| \Phi_F|_{H^1(\F_F^{\eps})} \leqslant C \eps^{- d / 2}$, we get with \eqref{I3}:
\begin{equation}
	\label{zbound} |z|_{H^1 (\Omega^{\eps})} \leqslant C \left( \eps^{- d} \sum_{F \in
		\Ein} J_F^2 \right)^{1 / 2} \leqslant C \eps^{5 / 2} \left( \sum_{F \in \Ein} \|
	\nabla_\tau I^0_{\eps}  \bar{u} \|^2_F \right)^{1 / 2} \leqslant C \eps^{3 / 2} \|
	\bar{u} \|_{H^1 (\Omega)} .
\end{equation}

\paragraph{Conclusion.} With $z$ chosen above, (\ref{zqident}) is satisfied,
so that (\ref{term3}) and (\ref{III1p12}) give $- III = III_1' + III_2$. Thus,
$| III |$ is bounded by the sum of (\ref{boundIII2}) and (\ref{boundIII1}).
Combining this with \ \eqref{boundI}, \eqref{zbound}, and inserting into \eqref{erreq} gives
\[ \tn{(u^{\eps} - \hat{u}^{\ast} + z, p^{\eps} - p^{\eps, 1})} \leqslant C \eps^{3/2} \| \bar{u} \|_{H^1 (\Omega)} . \]
This entails $|u^{\eps} - \hat{u}^{\ast} |_{H^1(\Omega^{\eps})}
\leqslant C \eps^{3 / 2} \| \bar{u} \|_{H^1 (\Omega)} + |z|_{H^1(\Omega^{\eps})} \leqslant C \eps^{3 / 2} \|
\bar{u} \|_{H^1 (\Omega)}$, thanks again to \eqref{zbound}. Lemma~\ref{lemPoin} gives $\|u^{\eps} -
\hat{u}^{\ast} \|_{\Omega^{\eps}} \leqslant C \eps^{5 / 2} \| \bar{u} \|_{H^1 (\Omega)}$. Since
$p^{\eps} - p^{\eps, 1} \in L^2_0 (\Omega^{\eps})$, Lemma~\ref{lemDiscPoin}
gives $$\|p^{\eps} - p^{\eps, 1} \|_{\Omega^{\eps}} \leqslant C \eps^{- 1} \tn{p^{\eps} - p^{\eps,
		1}} \leqslant C \sqrt{\eps} \, \| \bar{u} \|_{H^1 (\Omega)}\,.$$ Finally,
by \eqref{defansatz} and \eqref{I0},
\[ \|p^{\eps, 1} - p^{\ast} \|_{\Omega^{\eps}} \leqslant 2 \| \eps ({\pi}_i)_{\eps} I^0_{\eps}
\bar{u}_i \|_{\Omega^{\eps}} \leqslant C \eps \|I^0_{\eps}  \bar{u} \|_{\Omega} \leqslant C \eps \|
\bar{u} \|_{\Omega} . \]
This proves \eqref{mainest}.  \hfill$\Box$

\paragraph*{Use of large language models.} Claude (model Opus 5.5) was used to develop the central arguments in the proof of Theorem \ref{thmMain}. The article grew out of a first draft from July 2024, and has greatly benefited from interactions with Claude AI in many aspects. The most outstanding contributions were:
\begin{description}
\item{From human part:} The idea of using the finite element framework on the $\eps$-grid and the mesh-dependent norm (\ref{deftnqOld}) on the pressure. As mentioned in Remark \ref{remtnqOld}, this choice lead to technical complications and necessitated stronger regularity assumptions.  	
\item{From AI part:} Introducing the additional velocity correction $z$, as in (\ref{defz}).
\item{From human part again:} The idea of introducing gradient-like terms in the pressure norm (a priori inconsistent with the $L^2$ framework for the pressure). The norm (\ref{deftnqOld}) becomes (\ref{deftnq}). 
\end{description} 
Claude was also used for several other refinements, proof checks, bibliographical search, and language editing. The author has checked all mathematical arguments
and takes full responsibility for the content of the paper.

\paragraph*{Acknowledgements} The author warmly thanks Profs. Grégoire Allaire and François Murat for their interest in the first version of this article and enlightening  discussions. This work is part of the Franco-German project "Numerical homogenization techniques for multiscale Stokes problems" (PHC PROCOPE). The exchanges with German partners, in particular Moritz Hauck, were very helpful.   

\bibliographystyle{siam}
\bibliography{homog}
\end{document}